\documentclass[10pt]{article}

\usepackage{graphicx}
\usepackage{amsmath}
\usepackage{amsfonts}
\usepackage{amsthm}
\usepackage{amssymb}
\usepackage{url}
\usepackage{hyperref}  
\usepackage{verbatim}
\usepackage{color}

\newcommand{\figdir}{}

\def\mtx#1{\begin{bmatrix} #1 \end{bmatrix}}

\newcommand{\nul}{\operatorname{null}}

\def\ord#1{| #1 |}

\newcommand{\R}{\mathbb{R}}

\newcommand{\C}{\mathbb{C}}

\newcommand{\sRn}{S_n}
\newcommand{\Hn}{H_n}

\newcommand{\bx}{{\bf x}}

\newcommand{\G}{\mathcal{G}}
\newcommand{\SG}{\mathcal{S}(G)}
\newcommand{\SGp}{\mathcal{S}_+(G)}

\newcommand{\symp}{\mathcal{S}_+}
\newcommand{\hermp}{{\cal H}_+}
\newcommand{\rank}{\operatorname{rank}}
\newcommand{\mr}{\operatorname{mr}}
\newcommand{\hmrp}{\operatorname{mr}_+^\C}
\newcommand{\mrp}{\operatorname{mr}_+^\R}
\newcommand{\M}{\operatorname{M}}
\newcommand{\Mp}{\operatorname{M}_+^\R}
\newcommand{\hMp}{\operatorname{M}_+^\C}
\newcommand{\ZFS}{\operatorname{Z}}
\newcommand{\Zp}{\operatorname{Z}_+}
\newcommand{\OS}{\operatorname{OS}}
\newcommand{\Bmt}{B_m^t}
\newcommand{\HTH}{\overline{H_4(3)}}
\newcommand{\cc}{\operatorname{cc}}
\newcommand{\PC}{\operatorname{P}}

\newcommand{\cp}{\mathbin{\scriptscriptstyle\square}}

\newtheorem{thm}{Theorem}[section]
\newtheorem{prop}[thm]{Proposition}
\newtheorem{cor}[thm]{Corollary}
\newtheorem{lem}[thm]{Lemma}

\newtheorem{obs}[thm]{Observation}

\newtheorem{defn}[thm]{Definition}

\newtheorem{ex}[thm]{Example}
\newtheorem{quest}[thm]{Question}

\newcommand{\n}{\{1,\dots,n \}}
\newcommand{\x}{\times}

\newcommand{\bit}{\begin{itemize}}
\newcommand{\eit}{\end{itemize}}
\newcommand{\ben}{\begin{enumerate}}
\newcommand{\een}{\end{enumerate}}
\newcommand{\beq}{\begin{equation}}
\newcommand{\eeq}{\end{equation}}
\newcommand{\bea}{\begin{eqnarray*}}
\newcommand{\eea}{\end{eqnarray*}}
\newcommand{\bpf}{\begin{proof}}
\newcommand{\epf}{\end{proof}\ms}
\newcommand{\ms}{\medskip}

\title{Zero forcing parameters and minimum rank problems
\thanks{Part of this research was done at the American Institute of Mathematics SQuaRE,``Minimum Rank of Symmetric Matrices described by a Graph," and the authors thank AIM and NSF for their support.}}

\author{Francesco Barioli\thanks{Department of Mathematics, University of Tennessee at Chattanooga, Chattanooga TN, 37403
(francesco-barioli@utc.edu).} \and Wayne Barrett\thanks{Department of Mathematics, Brigham Young University, Provo UT 84602 (wayne@math.byu.edu).} \and Shaun M. Fallat\thanks{Department of Mathematics and Statistics,
University of Regina, Regina, SK, Canada
(sfallat@math.uregina.ca). Research supported in part by an
NSERC Discovery grant.}
\and \newline H. Tracy Hall\thanks{Department of Mathematics, Brigham Young University, Provo UT 84602 (H.Tracy@gmail.com).}
 \and Leslie Hogben\thanks{Department of Mathematics, Iowa State University, Ames, IA 50011, USA (lhogben@iastate.edu) and American Institute of Mathematics, 360 Portage Ave,
Palo Alto, CA 94306 (hogben@aimath.org).}   \and  Bryan  Shader\thanks{Department
of Mathematics, University of Wyoming,  Laramie, WY 82071, USA
(bshader@uwyo.edu).} \and \newline P. van den Driessche\thanks{Department of Mathematics and Statistics, University of Victoria, Victoria, BC, V8W 3R4, Canada
 (pvdd@math.uvic.ca).  Research supported in part by an
NSERC Discovery grant.} \and Hein van der Holst\thanks{School of Mathematics, Georgia Institute of Technology, Atlanta, GA 30332-0160, USA (holst@math.gatech.edu). On leave from Eindhoven University of Technology.}
 }

\date{}

\begin{document}
\maketitle

\begin{abstract}
The zero forcing number $\ZFS(G)$,
which is the minimum number of vertices in a zero
forcing set of a graph $G$,  is used to study
the maximum nullity/minimum rank of the family
of symmetric matrices described by $G$.  It is
shown that for a connected graph of order at
least two, no vertex is in every zero forcing set.
The positive semidefinite zero forcing number
$\Zp(G)$ is introduced, and shown to be equal to
$\ord G-\OS(G)$, where $\OS(G)$ is the recently defined ordered set number that is a lower bound for minimum positive semidefinite rank.  The positive semidefinite zero forcing
number is applied to the computation of  positive
semidefinite minimum rank of certain graphs.  An example of a graph for which the real positive symmetric semidefinite minimum rank is greater than the  complex Hemitian positive semidefinite minimum rank is presented.
\end{abstract}


\section{Introduction}\label{sintro}

 The {\em minimum rank problem} for a (simple) graph asks for the determination of  the minimum rank among all real symmetric matrices with the zero-nonzero pattern of off-diagonal entries  described by  a given  graph  (the diagonal of the matrix is free); the maximum nullity of the graph is the maximum nullity over the same set of matrices. This problem  arose from the study of possible eigenvalues of real symmetric matrices described by a graph and has received considerable attention over the last ten years (see \cite{FH} and references therein).  There has also been considerable interest in the related {\em positive semidefinite minimum rank problem}, where the minimum rank is taken over (real or complex Hermitian) positive semidefinite matrices described by a graph (see, for example, \cite{Netal08chord, CdVnu, Netal09OS, vdH03, vdH09, MNZ}).

 Zero forcing sets and the zero forcing number  were introduced in \cite{AIMgroup}.  The zero forcing number is a useful tool for determining the minimum rank of structured families of graphs and small graphs, and is motivated by simple observations about null vectors of
 matrices.  The zero forcing process is the same as graph infection used by physicists to study control of quantum systems \cite{graphinfect}, and the zero forcing number is becoming a graph parameter of interest in its own right.

 A {\em graph} $G=(V_G,E_G)$ means a simple undirected graph (no loops, no multiple edges) with a finite nonempty set of vertices $V_G$ and edge set $E_G$ (an edge is a two-element subset of vertices).   All matrices discussed are Hermitian; the set of  real symmetric $n\times n$ matrices is denoted by $\sRn$ and the set of  (possibly complex) Hermitian $n\times n$ matrices is denoted by $\Hn$.
For $A\in \Hn$,
the {\em graph} of $A$, denoted by $\G(A)$, is the graph with vertices
$\{1,\dots,n \}$  and edges $\{ \{i,j \} : ~a_{ij} \ne 0,  1 \le i <j \le n \}$.  Note that the diagonal of $A$ is ignored in determining $\G(A)$.  The study of minimum rank has focused on real symmetric matrices (or in some cases, symmetric matrices over a field other than the real numbers), whereas much of the work on positive semidefinite minimum rank  involves (possibly complex) Hermitian matrices.  Whereas it is well known that using complex Hermitian matrices can result in a lower minimum rank than using real symmetric matrices, one of the issues in the study of minimum positive semidefinite rank has been whether or not using only real matrices or allowing complex matrices matters to minimum positive semidefinite rank. Example \ref{4hub3ex} below shows that complex Hermitian positive semidefinite minimum rank can be strictly lower than real symmetric positive semidefinite minimum rank.

Let $G$ be a graph. The {\em set of real symmetric matrices described by  $G$} is
\[  \SG=\{A\in\sRn : \G(A)=G\}. \]
The {\em minimum rank} of  $G$  is  \[\mr(G)=\min\{ \rank A : A\in\SG\}\]
and the {\em maximum nullity of $G$} is
\[\M(G)=\max\{\nul A: A\in \SG\}.\]
 Clearly
$\mr(G)+\M(G)=\ord G,$
where the {\em order} $\ord G$ is the number of vertices of $G$.   The {\em set of real positive semidefinite matrices described by  $G$} and the {\em set of Hermitian positive semidefinite matrices described by  $G$} are, respectively,
\bea \SGp & =& \{A\in\sRn : \G(A)=G \mbox{ and $A$ is positive semidefinite}\}\\
\hermp(G)&=&\{A\in H_n : \G(A)=G \mbox{ and  $A$ is positive semidefinite}\}. \eea
The {\em minimum positive semidefinite  rank}  of $G$ and {\em minimum Hermitian positive semidefinite  rank}  of $G$ are, respectively,
\[\mrp(G)=\min\{ \rank A : A\in\SGp\} \mbox{ and } \hmrp(G)=\min\{ \rank A : A\in\hermp(G)\}.\]
The {\em maximum positive semidefinite nullity of $G$} and the {\em maximum Hermitian positive semidefinite nullity of $G$} are, respectively,
\[\Mp(G)=\max\{\nul A: A\in \SGp\} \mbox{ and } \hMp(G)=\max\{\nul A: A\in \hermp(G)\}.\]
Clearly
$\mrp(G)+\Mp(G)=\ord G$ and $\hmrp(G)+\hMp(G)=\ord G.$
There are a variety of symbols in the literature (see, for example,
\cite{Netal08chord, MNZ}) for these
parameters, including ${\rm msr}(G)$ and ${\rm hmr}_+(G)$ for
what we denote by $\hmrp(G)$.    Clearly $\Mp(G)\le\M(G)$ and $\mr(G)\le\mrp(G)$ for every graph $G$, and it is well known that these inequalities can be strict (for example, any tree $T$ that is not a path has $\mr(T)<\mrp(T)$).

We need some additional  graph terminology.  The {\em complement} of a graph $G=(V,E)$ is the
graph $\overline{G}=(V,\overline{E})$, where $\overline{E}$
consists of all two element sets from $V$ that are not in $E$.
We denote  the complete graph on $n$ vertices by $K_n$; a complete graph is also called a clique.
The {\em degree} of vertex $v$ in graph $G$ is the number of edges incident with  $v$, and the minimum degree of the vertices of  $G$ is denoted by $\delta(G)$.  A set of subgraphs of $G$, each of which is a clique and such that every edge of $G$ is
contained in at least one of these cliques, is called a {\em clique covering} of $G$.
 The {\em clique covering number} of $G$, denoted by $\cc(G)$, is the smallest number of
cliques in a clique  covering of $G$.

\begin{obs}\label{mrpcc}  {\rm \cite{FH}} For every graph $G$, $\mrp(G)\le\cc(G)$, so $\ord G -\cc(G)\le \Mp(G)$.
\end{obs}

For an $n\x n$ matrix $A$ and $W \subseteq \n$, the principal submatrix $A[W]$ is the submatrix of
$A$ lying in the rows and columns that have indices in $W$.   For a graph $G=(V_G,E_G)$ and $W\subseteq V_G$, the {\em induced subgraph} $G[W]$  is the graph with vertex set $W$ and edge set $\{\{v,w\} \in E_G : v,w \in W\}$.     The induced subgraph $\G(A)[W]$ of the graph of $A$ is naturally associated with the graph of the the principal submatrix for $W$, i.e., $\G(A[W])$.
The subgraph induced by $\overline W =V_G \setminus W$ is usually denoted by $G-W$,
 or in the case $W$ is a singleton $\{v\}$, by $G-v$.

The {\em path cover number} $\PC(G)$ of $G$ is the smallest positive integer $m$ such that there are $m$ vertex-disjoint induced paths  $P_1,\dots, P_m$  in $G$ that cover all the vertices of $G$ (i.e.,   $V_G=\dot{\cup}_{i=1}^m V_{P_i}$). A graph is {\em planar} if it can be drawn in
the plane without crossing edges.  A graph is {\em outerplanar} if it has  such a drawing with a face
that contains all vertices.  Given two graphs $G$ and $H$, the {\em Cartesian product}
of $G$ and $H$, denoted $G \cp H$, is the graph whose vertex set is the
Cartesian product of $V_G$ and $V_H$, with an edge
between two vertices exactly when they are identical
in one coordinate and adjacent in the other.

 Let $G=(V_G,E_G)$ be a  graph.
 A subset $Z\subseteq V_G$ defines an initial set of  black vertices (with all the vertices not in $Z$ white), called a {\em coloring}.
 There are no constraints on permissible colorings; instead there are constraints on how new colorings can be derived.
 The {\em  color change rule} (for the zero forcing number) is to change the color of a white vertex $w$ to black if $w$ is the unique white neighbor of a black vertex $u$; in this case we say $u$ {\em forces} $w$  and write $u\to w$.  Given a coloring of  $G$, the  {\em  derived set}  is the set of black vertices obtained by  applying the color change rule until no more changes are possible.
 A  {\em  zero forcing set}   for  $G$  is a subset of vertices $Z$ such that if initially the vertices in $Z$ are colored black and the remaining vertices are colored white, the derived set is $V_G$. The  {\em zero forcing number} $\ZFS(G)$ is the minimum of $\ord Z$ over all   zero forcing sets $Z\subseteq V_G$.

\begin{thm}\label{propZ}{\rm \cite[Proposition 2.4]{AIMgroup}} For any graph $G$, $\M(G) \leq\ZFS(G)$.  \end{thm}

Suppose $S=(v_1, v_2, \ldots, v_m)$ is an ordered subset of vertices from a
given graph $G$. For each $k$ with $1 \leq k \leq m$,
let $G_{k}$ be the subgraph of $G$ induced by
$\{v_1, v_2, \ldots, v_k\}$, and let
$H_k$ be the connected component of $G_k$ that contains $v_k$.
If for each $k$, there exists a vertex $w_k$ that satisfies:
$w_k \neq v_l$ for $l \leq k$, $\{w_k, v_k \} \in E$, and
$\{w_k, v_s \} \not\in E$, for all $v_s$ in $H_k$ with $s \neq k$,
then $S$ is called an
{\em ordered set of vertices} in $G$, or an OS-set. As defined in \cite{Netal09OS}, the {\em OS number}
of a graph $G$, denoted by $\OS(G)$, is the maximum of  $|S|$ over all OS-sets $S$ of $G$.

\begin{thm}\label{OSmr}{\rm \cite[Proposition 3.3]{Netal09OS}} For any  graph $G$, $\OS(G) \leq\hmrp(G)$.  \end{thm}

In Section \ref{sZFN} we establish several properties of the zero forcing number, including the nonuniqueness of zero forcing sets.  In Section \ref{sPSDZ} we introduce the positive semidefinite zero forcing number as an upper bound for maximum positive semidefinite nullity,  show that  the sum of the positive semidefinite zero forcing number and the OS number is the order of the graph, and apply the positive semidefinite zero forcing number to the computation of  positive semidefinite minimum rank.  Section \ref{smsr} provides the first example showing that $\mrp(G)$ and $\hmrp(G)$ need not be the same (described as unknown in \cite{FH}).


\section{Properties of the zero forcing number}\label{sZFN}

In this section, we establish several properties of  the zero forcing number, including the non-uniqueness of zero forcing sets and its relationship to path cover number.  
 We need some additional definitions related to the zero forcing number.
 \begin{defn} {\rm A {\em minimum zero forcing set} is a zero forcing set $Z$ such that $\ord Z = \ZFS(G)$.
}\end{defn}

Zero forcing chains of digraphs were defined in  \cite{square1}.  We give an analogous definition for graphs. 
\begin{defn} {\rm Let $Z$ be a  zero forcing set  of a graph $G$. \bit
\item Construct the  derived set, recording the forces in the order in which they are performed.
This is the {\em chronological list of forces}.
\item
A {\em forcing chain} (for this particular chronological list of
forces) is a sequence of vertices $(v_1,v_2,\dots,v_k)$ such that
for $i=1,\dots,k-1$,  $v_i \to v_{i+1}$. \item A {\em maximal
forcing chain} is a forcing chain  that is not a proper subsequence
of another zero forcing chain. \eit Note that a zero forcing chain
can consist of a single vertex $(v_1)$,  and such a chain is maximal
if $v_1\in Z$ and $v_1$ does not perform a force. }\end{defn}

As noted in \cite{AIMgroup}, the  derived set of a given set of
black vertices is unique; however, a chronological list of forces (of one particular zero forcing set) usually is not.
At Rocky Mountain Discrete Mathematics Days held Sept. 12 -- 13, 2008 at the University of Wyoming, the following questions were raised.

\begin{quest}\label{q1}  Is there a graph  that has a unique minimum zero forcing set?
\end{quest}
\begin{quest}\label{q2} Is there a graph $G$ and a vertex $v\in V_G$ such  that $v$ is in every minimum zero forcing set?
\end{quest}

We  show the answers to both these questions are negative for nontrivial connected graphs.

\begin{defn} {\rm Let $Z$ be a  zero forcing set of a graph $G$.
 A {\em reversal}  of $Z$ is the set of last vertices of the maximal zero forcing chains of a chronological list of forces. }\end{defn}

Each vertex can force at most one other vertex and can be forced by at most one other vertex, so the maximal forcing chains are disjoint, and the elements of $Z$ are the initial vertices of the maximal forcing  chains.  Thus the cardinality of a reversal of $Z$ is the same as the cardinality of $Z$.

\begin{thm}  If $Z$ is  a zero forcing set of $G$ then so is any reversal of $Z$.
\end{thm}
\bpf Write the chronological list of forces in reverse order, reversing each force (call this the reverse chronological list of forces) and let the reversal of $Z$ for this list be denoted $W$.  We show the reverse chronological list of forces is a valid list of forces for $W$. Consider the first ``force" $u\to v$ on the reverse chronological list. All neighbors of $u$ except $v$ must be in $W$, since when the last force $v\to u$ of $Z$ was done, each of them had the white neighbor $u$ and thus did not  force any vertex previously (in the original chronological list of forces).  Thus $u\to v$ is  a valid force for $W$.  Continue in this manner or use induction on $\ord G$.
\epf

\begin{cor}  No connected graph of order greater than one has a unique minimum zero forcing set.
\end{cor}

\begin{lem}\label{noallblack}  Let $G$ be a connected graph of order greater than one and let $Z$ be a minimum zero forcing set.  Every $z\in Z$ has a neighbor $w\not\in Z$.
\end{lem}
\bpf Suppose not.  Then there is a vertex $z\in Z$ such that every neighbor  of $z$ is in $Z$ (and $z$ does have at least one neighbor $v$).   Since $z$ cannot perform a force, $z$ is in the reversal $W$ of $Z$.  Using the reversed maximal forcing chains, no neighbor of $z$ performs a force. So $W\setminus\{z\}$ is a zero forcing set of smaller cardinality, because after every vertex except $z$ is black, $v$ can force $z$. \epf

\begin{thm}\label{noallZ}  If $G$ is a connected graph of order greater than one, then \[{\bigcap}_{Z \in ZFS(G)} Z = \emptyset,\]  where  $ZFS(G)$ is the set of all minimum zero forcing sets of $G$.
\end{thm}
\bpf Suppose not.  Then there exists $v\in {\cap}_{Z \in ZFS(G)} Z $.  In particular, for each $Z$ and each reversal $W$ of $Z$, $v$ is in both $Z$ and $W$.  This means that there is a maximal forcing chain consisting of only $v$, or in other words $v$ does not force any other vertex.  

Let $Z$ be a zero forcing set.  If there is no chronological list of forces in which  a neighbor of $v$ performs a force, then replace $Z$ by its reversal (since, by Lemma \ref{noallblack}, $v$ originally had a white neighbor $u$, in the reversal $u$ performs a force).  Let $u\to w$ be the first force in which the forcing vertex $u$ is a neighbor of $v$.  We claim that $Z\setminus \{v\}\cup \{w\} $ is a zero forcing set for $G$.  The  forces can proceed until $u$ is encountered as a forcing vertex.  At that time, replace $u\to w$ by $u\to v$, and then continue as in the original chronological list of forces.
\epf

  Next we show that for any graph the zero forcing number is an upper bound for the path cover number.

\begin{prop}\label{propZP} For any graph $G$, $\PC(G)\le \ZFS(G)$.
\end{prop}
\bpf Let $Z$ be a zero forcing set. The vertices in a  forcing chain induce a path in $G$ because  the forces in a forcing chain   occur chronologically  in the order of the chain (since only a black vertex can force). The maximal forcing chains are disjoint, contain all the vertices of $G$, and the elements of the set $Z$  are the initial vertices of the maximal forcing  chains.  Thus $\PC(G)\le \ord Z$.  By choosing a minimum zero forcing set $Z$, $\PC(G)\le \ZFS(G)$.
\epf

 In \cite{JLD} it was shown that for a tree $T$, $\PC(T)=\M(T),$ and
 in \cite{AIMgroup} it was shown that for a tree, $\PC(T)=\ZFS(T)$
 (and thus $\M(T)=\ZFS(T)$).
In \cite{BFH}  it was shown that for graphs in general, $\PC(G)$ and $\M(G)$
are not comparable.
However,  Sinkovic has established the following relationship for outerplanar graphs:
    If  $G$ is an  outerplanar graph, then     $  \M(G)\le \PC(G)$ \cite{Sink}.  The next example shows that neither outerplanar graphs nor 2-trees require $\M(G)=\ZFS(G)$ or $\PC(G)=\ZFS(G)$ (a 2-tree  is constructed inductively by starting with a
 $K_3$ and connecting each new vertex to 2 adjacent existing
  vertices).

\begin{ex}\label{pinwheelex}{\rm Let $G_{12}$ be the graph  shown in Figure \ref{fig:pinwheel}, called the  pinwheel on 12 vertices. Note that $G_{12}$ is an outerplanar 2-tree.
 \begin{figure}[!h]\begin{center}
\scalebox{.6}{\includegraphics{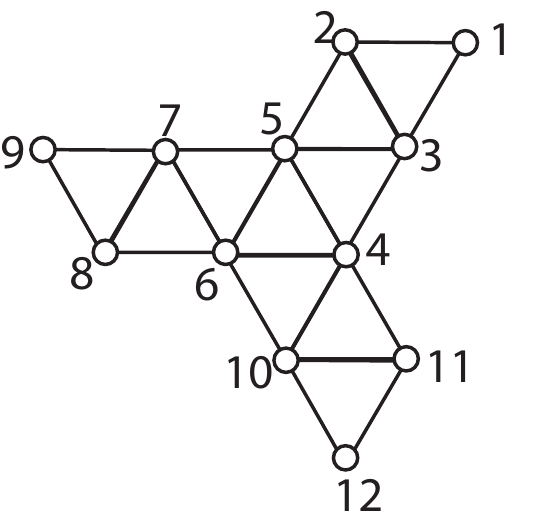}}
\caption{The graph $G_{12}$ for Example \ref{pinwheelex}, the pinwheel on $12$ vertices}\label{fig:pinwheel}
\end{center}
\end{figure}
The set  $\{1,2,6,10\}$ is a  zero forcing set for $G_{12}$, so  $\ZFS(G_{12})\le 4$.  We show that $\ZFS(G_{12})\ge 4$, which implies $\ZFS(G_{12})=4$.
  Suppose to the contrary that  $Z$ is a  zero forcing
  set for $G_{12}$ and $\ord Z=3$.  To start the
  forcing, at least two of the vertices must be
  in one of the sets $\{1,2,3\}, \{7,8,9\}, \{10,11,12\}$;
  without loss of generality, assume that
  two or three black vertices are in $\{1,2,3\}$.  Then after several forces the vertices $\{1,2,3, 4, 5\}$ are black, and at most one additional vertex $v\not\in\{1,2,3,4,5\}    $ is in $Z$.  To perform another force with only one more black vertex $v$, either $6$ or $7$ must be black, and $5$ can force the other, but then no additional forces can be performed, so $Z$ was not a  zero forcing set for $G_{12}$.
      Clearly  $G_{12}$ can be covered by 9 triangles, so
      $\cc(G_{12})\le 9$ and $\Mp(G_{12})\ge 3$, by Observation \ref{mrpcc}.
 It is easy to find a path covering of 3 paths,
 so  $\Mp(G_{12})=\M(G_{12})=\PC(G_{12})=3$ and $\mrp(G_{12})=\mr(G_{12})=\cc(G_{12})=9$. 
 Since $G_{12}$ is chordal,
 $\hmrp(G_{12})=\cc(G_{12})$ \cite{Netal08chord},
 and thus $\hMp(G_{12})=3$.}\end{ex}


\section{The positive semidefinite zero forcing number}\label{sPSDZ}

  In this section, we introduce the positive definite zero forcing number, relate it to maximum positive semidefinite nullity and to the OS number, and apply it to compute  maximum positive semidefinite nullity of several families of graphs.
  The definitions and terminology for zero forcing (coloring, derived set, etc.) are the same as for the zero forcing number $\ZFS(G)$, but the color change rule is different.

\begin{defn}\label{defZp}$\null$ {\rm  \bit
\item   The {\em positive semidefinite color change rule} is:\\
 Let $B$ be the set consisting of all the black vertices.  Let $W_1,\dots, W_k$ be the sets of vertices of the $k$ components of $G-B$ (note that it is possible that $k=1$).  Let $w\in W_i$.
If $u\in B$ and $w$ is the only white neighbor of  $u$ in $G[W_i\cup B]$, then change the color of $w$ to black.\item The  {\em  positive semidefinite zero forcing number of a  graph $G$}, denoted by $\Zp(G)$,   is the minimum of $\ord X$ over all positive semidefinite   zero forcing sets $X\subseteq V_G$ (using the positive semidefinite color change rule).  \eit }\end{defn}

Forcing using the positive semidefinite color change rule can be thought of as decomposing the graph into a union of certain  induced subgraphs and using ordinary zero forcing on each of these induced subgraphs.
The application  of the positive semidefinite color change rule  is illustrated in the next example.

\begin{ex}\label{extree}  {\rm   Let $T$ be a tree.  Then $\Zp(T)=1$, because
 any one vertex $v$ is a positive semidefinite zero forcing set.
 Formally, this can be established by induction on $\ord T$:
 If $v$ is a leaf, it forces its neighbor; if not a
 decomposition takes place.  In either case smaller
 tree(s) are obtained.   It has been known for a long time (see, for example, \cite{FH}) that  $\hMp(T)=1$,
  but the use of $\Zp$ provides an easy proof of this result, because $\hMp(T)=1$ is an immediate consequence of $\Zp(T)=1$ by Theorem \ref{propZp} below.
 }\end{ex}

\begin{obs}  Since any zero forcing set is a positive definite zero forcing set, \end{obs}\[\Zp(G)\le \ZFS(G).\]

\begin{ex}\label{expinwheel2}  {\rm The pinwheel $G_{12}$ shown in
Figure \ref{fig:pinwheel} has $\Zp(G_{12})=3=\Mp(G_{12})$
because $X=\{4,5,6\}$  is a positive semidefinite
zero forcing set ($G_{12}-X$ is disconnected,
and $X$ is a zero forcing set for $G[\{1,2,3,4,5,6\}]$, etc.).}\end{ex}

For any graph $G$ that is the disjoint union of connected components
$G_i, i=1,2,\ldots, k$, $\Zp(G) = \sum_{i=1}^{k} \Zp(G_i)$ (the analogous
results for $\M, \Mp, \hMp$ and $\ZFS$ are all well known).

\begin{thm}\label{propZp} For any graph $G$, $\hMp(G) \leq\Zp(G)$.
\end{thm}
\bpf Let  $A\in\hermp(G)$ with $\nul A=\hMp(G)$. Let $\bx=[x_i]$ be a nonzero vector in  $\ker A$.  Define $B$ to be the set of indices $u$ such that $x_u=0$ and  let $W_1,\dots, W_k$ be the sets of vertices of the $k$ components of $G-B$.
We claim that   in  $G[B\cup W_i]$, $w\in W_i$ cannot be the unique
neighbor  of any vertex $u\in B$.  Once the claim is
established,  if  $X$ is a positive semidefinite zero forcing set for $G$,
then  the only  vector in $\ker A$
with zeros in positions indexed by $X$ is the zero vector, and
thus $ \hMp(G)\le\Zp(G)$.

To establish the claim, renumber the vertices so that the vertices of $B$ are last, the vertices of $W_1$ are first, followed by the vertices of $W_2$, etc.  Then $A$ has the block form
\[A=\mtx{
A_1 & 0 & \dots & 0 & C_1^*\\
0 & A_2 & \dots & 0 & C_2^*\\
\vdots & \vdots & \ddots & \vdots\\
0 & 0 & \dots &  A_k & C_k^*\\
C_1 & C_2 & \dots & C_k & D
}.\]
Partition $\bx$ conformally as $\bx=[\bx_1^T,\dots,\bx_k^T,0]^T$, and note that all entries of $\bx_i$ are nonzero, $i=1,\dots,k$.
Then $A\bx=0$ implies $A_i\bx_i=0, i=1,\dots,k$.  Since $A$ is positive semidefinite, each column in $C_i^*$ is in the span of the columns of $A_i$ by the column inclusion property of Hermitian positive semidefinite matrices \cite{FJ}.  That is, for $i=1,\dots k$, there exists $Y_i$ such that $C_i^*=A_iY_i$.  Thus  $ C_i \bx_i =Y_i^*A_i\bx_i =0$, and $w\in W_i$ cannot be the unique neighbor  in $W_i$ of any vertex $u\in B$.
\epf

Theorem \ref{propZp} is also a consequence of Theorem \ref{OSZpG}  below and Theorem \ref{OSmr} above, 
 but using that as a justification obscures  the motivation for the definition and the connection between zero forcing and null vectors that is given in the short direct proof.

In \cite[Theorem 2.10]{MNZ} it is shown that $\ord G -\ZFS(G)\le OS(G)$.  A similar method can be used to show an a more precise relationship between $\Zp$ and the OS number.

\begin{thm}\label{OSZpG} For any graph $G=(V,E)$ and any ordered set $S$, $V\setminus S$ is a positive semidefinite forcing set for $G$, and for any positive semidefinite forcing set $X$ for $G$, there is an order that makes $V\setminus X$ an ordered set for $G$.  Thus $\Zp(G) + OS(G)=|G| $.
\end{thm}
\bpf Let $X$ be a positive semidefinite zero forcing set for $G$  such that $|X| = \Zp(G)$. Let $v_{i}$ be the
vertex colored black by the $i$th application of the positive semidefinite color change rule. We  show that $S = (v_{t}, v_{t-1}, \ldots, v_{1})$ is an OS set for
$G$, where $t=|G|-\Zp(G)$. Further define $X_{0}=X$, and
$X_{i+1} = X_{i} \cup \{v_{i+1}\}$, for $i=0,1,\ldots, t-1$.
For each $v_{i}$, since it was initially white and then colored black on the
$i$th
application of the positive semidefinite color change rule, there exists a vertex $w_{i} \in X_{i}$
(the current black vertices) such that $v_{i}$ is the only neighbor in
the subgraph of $G$ induced by $X_{i} \cup H_{1}$, where
the subgraph $G \setminus X_{i}$ has components $H_{1}, H_{2}, \ldots H_{p}$
with $v_i \in H_1$. Since $X$ is a positive semidefinite zero forcing set, no
other vertex from the set $\{v_{i+1}, v_{i+2}, \ldots, v_{t}\}$ (the remaining
white vertices) can be in $H_{1}$ and be a neighbor of $w_i$. Hence the set
$(v_{t}, v_{t-1}, \ldots, v_{1})$ is an OS-set. Therefore $t \leq OS(G)$.
Thus
\beq |G| - \Zp(G)\leq OS(G).\label{OSZp1}\eeq

For the converse, we use the fact that if
$S=(v_1, v_2, \ldots v_m)$ is an OS set, then the set $S \setminus \{v_m\}$
is also an OS set.
Suppose $S=(v_1, v_2, \ldots v_m)$ is an OS set with $|S|=OS(G)$. Then
we claim that $V\setminus S$ is a positive semidefinite
zero forcing set.
So color the vertices $V\setminus S$ black, and suppose the subgraph $G_m$
induced by the vertices of $\{v_1,\dots,v_m\}$ 
has components induced by $U_1, U_2, \ldots, U_\ell$.
Let $v_m \in U_1$.
Since $S$ is an OS-set there exists a vertex
$w_m \in V \setminus S$ such that $w_m v_m \in E$ and $w_m v_s \not\in E$
for all other $v_s \in U_1$. This implies that $v_m$ can be colored black
under the positive semidefinite color change rule.
Since $S \setminus \{v_m\}$
 is also an OS-set for $G$,
we may continue this argument and deduce that $V \setminus S$ is a positive semidefinite zero
forcing set. Hence
\beq |G| - OS(G) = |V \setminus S| \geq Z_{+}(G),\label{OSZp2}\eeq
as the positive semidefinite zero forcing number  is
defined as a minimum over all such zero forcing sets.
From (\ref{OSZp1}) and (\ref{OSZp2}), $Z_{+}(G) + OS(G)= |G|$.
\epf

\begin{cor} For every graph $G$, \end{cor}\[ \delta(G) \leq \Zp(G). \]

\bpf  By \cite[Corollary 2.19]{MNZ}, $OS(G)\le \ord G - \delta(G)$.
Combining this with Theorem \ref{OSZpG} gives the result. \epf

Another consquence of Theorem \ref{OSZpG} is that there are
examples of graphs for which $\Zp$ may not be equal to $\hMp$. For
example, in \cite{MNZ} it was shown that the M\"obius Ladder on 8 vertices,
sometimes denoted by $ML_8$ or $V_8$, satisfies $OS(ML_8)=4$ and $\hmrp(ML_8) = 5$.
In this case, by Theorem \ref{OSZpG}, it follows  that $\Zp(ML_8)=4$, and
hence $\Zp(ML_8) > 3= \hMp(ML_8)$.

In \cite{AIMgroup}, the zero forcing number was used to establish the minimum rank/maximum nullity of numerous families of graphs.  The positive semidefinite zero forcing number is equally effective. Here we apply it to two families of graphs. 
The set of vertices associated with (the same) positive semidefinite zero forcing set in each copy of $G$ is a positive semidefinite zero
forcing set for $G\cp H$.

\begin{prop}\label{cpZpbd} For all graphs $G$ and $H$,
$\Zp(G\cp H)\le \min\{\Zp(G)\ord H, \Zp(H)\ord G\}$.
\end{prop}

\begin{cor}\label{cortree}  If $T$ is a tree and $G$ is a graph, then
$\Zp(T\cp G)\le \ord G$.
\end{cor}

\begin{thm}\label{treeKr}  If $T$ is a tree of order at least two,
then $\Mp(T\cp K_r)=\hMp(T\cp K_r)=\Zp(T\cp K_r)=r$.
\end{thm}
\bpf Let $T$ be a tree of order $n\ge 2$. By Corollary \ref{cortree}, $\Zp(T\cp K_r)\le r$.  We show $r\le \Mp(T\cp K_r)$ by constructing a matrix $A\in\symp(T\cp K_r)$ of rank at most $(n-1)r$, and the result then follows from Theorem  \ref{propZp}.
The construction is by induction on $n$.  Let $P_2$ denote the path on $2$ vertices.  To show that $\mrp(P_2\cp K_r)=r$, choose a nonsingular matrix $M\in\symp(K_r)$ such that $M^{-1}\in\symp(K_r)$ (for example, $M=I+J$, where $I$ is the identity matrix and $J$ is the all 1s matrix).  Then $B=\mtx{M & I\\I & M^{-1}}\in \symp(P_2\cp K_r)$ and $\rank B=\rank M=r$. Without loss of generality, in $T$ vertex $n$ is adjacent only to vertex $n-1$.  We order the vertices $(i,j)$ of $T\cp K_r$ lexicographically.  By the induction hypothesis, there is  a matrix $C\in\symp((T-n)\cp K_r)$ such that $\rank C = (n-2)r$; let $C'=C\oplus 0_{r\x r}$.  Using $B\in\symp(P_2\cp K_r)$ already constructed with rank $r$,  let $B'=0_{(n-2)r\x (n-2)r}\oplus B$.   Then for $\alpha\in\R$ chosen to avoid cancellation, $A=C'+\alpha B'\in\symp(T\cp K_r)$ and $\rank A\le (n-2)r+r=(n-1)r$.
\epf

A {\em book} with $m\ge 2$ pages, denoted $B_m$ \cite[p. 14]{Gal},
is  $m$ copies of a 4-cycle with one edge in common, or
equivalently, $B_m=K_{1,m}\cp P_2$, where $K_{1,m}$ is the complete
bipartite graph with partite sets of 1 and $m$ vertices. For $m \ge
2$, $t \ge 3$, we call $m$ copies of a $t$-cycle with one edge in
common a \textit{generalized book}, denoted by $B_m^t$
(obviously, $B_m=B_m^4$).

\begin{prop}\label{book}  If $\Bmt$ is a generalized book, then
$\Mp(\Bmt)=\hMp(\Bmt)=\Zp(\Bmt)=2$.
\end{prop}
\bpf  The two vertices in the common edge are a
positive semidefinite zero forcing set, so $\Zp(\Bmt)\le 2$.  Thus by Theorem \ref{propZp}, $\hMp(B_m^t) \le 2$.  Since $\Bmt$ is not a tree, $\Mp(\Bmt)\ge 2$ \cite{vdH03}.
\epf


\section{Real versus complex minimum positive semidefinite rank}\label{smsr}  Clearly  $\hmrp(G)\le\mrp(G)$ for every graph $G$. Previously it was not known whether $\hmrp(G)$ could differ from $\mrp(G)$ \cite[p. 578]{FH}.
In this final section we provide an example of a graph for which these parameters are not identical.

\begin{ex}\label{4hub3ex}{\rm The  ``$k$-wheel with 4 hubs" (for $k$ at least 3) is the graph on $4k+4$ vertices such that the outer cycle has $4k $ vertices, and each of the 4 hubs is attached to every 4th vertex of the cycle, and no others; this graph is denoted $H_4(k)$, and $H_4(3)$ is shown in Figure \ref{fig:4hub3}. This family arose in Hall's investigation of graphs having minimum rank 3 \cite{Hall4hub}. We show $\hmrp(\HTH)=3$ and $\mrp(\HTH)=4$.
 \begin{figure}[!h]\begin{center}
\scalebox{.4}{\includegraphics{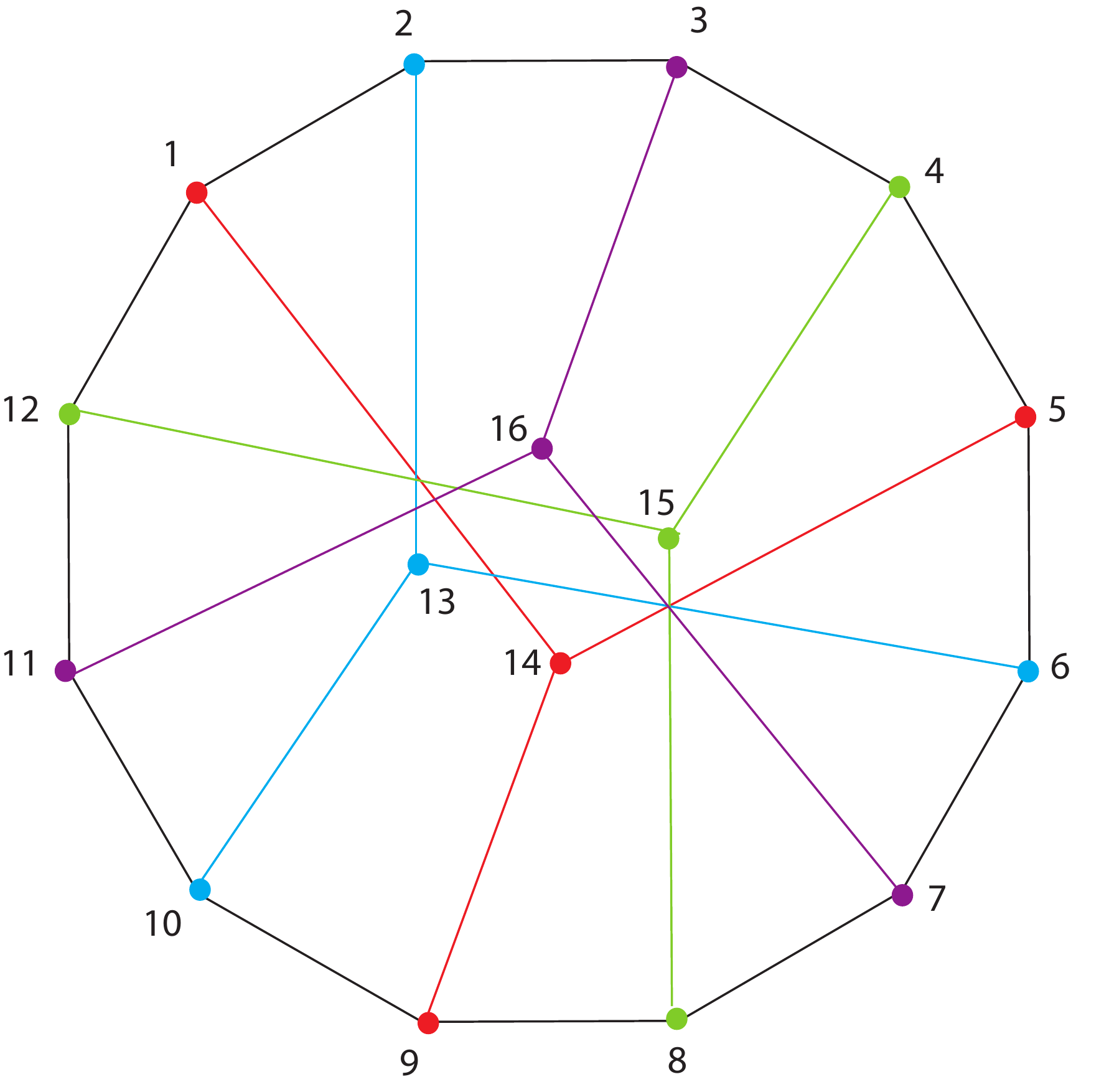}}
\caption{The  the  3-wheel on 4 hubs, $H_4(3)$, for Example \ref{4hub3ex}}\label{fig:4hub3}
\end{center}
\end{figure}
As numbered in Figure \ref{fig:4hub3}, $H_4(3)$ is bipartite with partite sets consisting of the odd vertices and the even vertices.  By \cite[Theorem 3.1]{square1}, $\mrp(\HTH)=\mr^\R(Y_{\HTH})$ where $Y_{\HTH}$
 is the biadjacency zero-nonzero pattern of $\HTH$ and $\mr^\R(Y_{\HTH})$ is the asymmetric minimum rank over the real numbers (Theorem 3.1   applies to $\HTH$ because $H_4(3)$ is a  bipartite graph).
 The same method used to prove Theorem 3.1 also shows that  $\hmrp(\HTH)=\mr^\C(Y_{\HTH})$ where $\mr^\C(Y_{\HTH})$ is the asymmetric minimum rank over the complex numbers (in \cite[Remark 3.2]{square1} it is noted that the method in Theorem 3.1 is valid for constructing a symmetric matrix over an infinite field, and the same reasoning applies to constructing a Hermian matrix over $\C$ by using Hermitian adjoints in place of transposes).   After scaling rows and columns, a minimum rank matrix having zero-nonzero pattern $Y_{\HTH}$ has the form
 \[A=\mtx{0 & 1 & 1 & 1 & 1 & 0 & 0 & 1 \\
 0 & 0 & 1 & a_{3,8} & a_{3,10} & a_{3,12} & a_{3,14} & 0 \\
 1 & 0 & 0 & a_{5,8} & a_{5,10} & a_{5,12} & 0 & a_{5,16} \\
 1 & a_{7,4} & 0 & 0 & a_{7,10} & a_{7,12} & a_{7,14} & 0 \\
 1 & a_{9,4} & a_{9,6} & 0 & 0 & a_{9,12} & 0 & a_{9,16} \\
 1 & a_{11,4} & a_{11,6} & a_{11,8} & 0 & 0 & a_{11,14} & 0 \\
 0 & 1 & 0 & a_{13,8} & 0 & a_{13,12} & a_{13,14} & a_{13,16} \\
 1 & 0 & a_{15,6} & 0 & a_{15,10} & 0 & a_{15,14} & a_{15,16}}.\]
where the displayed entries $a_{ij}$ are nonzero (real or complex) numbers. Since the principal submatrix in the first three rows and columns is nonsingular, $\rank A=3$ implies that rows 4 through 8 are linear combinations of rows 1 through 3. Computations show that the following assignments of variables are necessary:
\bea a_{5,8} & = &  \left(a_{3,8}-1\right) a_{7,4},\,\,
a_{5,10} = \left(a_{3,10}-1\right) a_{7,4}+a_{7,10},\,\,
   a_{5,12}  =  a_{3,12} a_{7,4}+a_{7,12},\,\,
     a_{5,16}  =  - a_{7,4},\\
  a_{7,14} & = &-   a_{3,14}  a_{7,4},\,\,
  a_{ 9, 16} = a_{9, 4} - a_{7, 4}, \,\,
  a_{9, 6} = a_{9, 4} , \,\,
  a_{9, 12} =   a_{3, 12} a_{7, 4} + a_{7, 12} -
   a_{3, 12} a_{9, 4} +   a_{3, 12} a_{9, 6}, \\
   a_{7, 10} &=&  a_{7, 4} - a_{3, 10} a_{7, 4} -  a_{9, 4}, \,\,
   a_{9, 4} = (1 - a_{3, 8}) a_{7, 4}, \,\,
   a_{11, 4} = a_{7, 4}, \,\,
a_{11, 14} =   a_{3, 14} ( a_{11, 6} - a_{11, 4}), \\
a_{7, 12} &=& -a_{3, 12} a_{11, 6}, \,\,
a_{11, 8} = a_{3, 8} a_{11, 6}, \,\,
a_{3, 8} =   a_{3, 10} (a_{7, 4} - a_{11, 6})/a_{7, 4}, \,\,
a_{13, 16} = 1, \,\,
a_{13, 14} = -a_{3, 14}, \\
a_{13, 12} &=& -a_{3,  12}, \,\,
a_{3, 10}= 1, \,\,
a_{13, 8} = a_{11, 6}/a_{7, 4}, \,\,
a_{15, 16} = -a_{7, 4}, \,\,
a_{15, 14} =  a_{3, 14} a_{15,   6}, \\
a_{15, 10} &=& -a_{11, 6} + a_{15,   6}, \,\,
  a_{11, 6} = a_{7, 4} + a_{15, 6}.
\eea
After making these assignments, rows 4 - 7 are linear combinations of rows 1, 2, and 3, and in order for row 8   to be a linear combinations of rows 1, 2, and 3, it is necessary and sufficient that
\beq 1+\frac{  a_{7,4}}{a_{15,6}}  + \left(\frac{  a_{7,4}}{a_{15,6}}\right)^2 = 0.\label{stickyeq}\eeq
Clearly (\ref{stickyeq}) has a solution if and only if the field contains a primitive third root of unity.  Thus
$\mr^\C(Y_{\HTH}) = 3$ whereas $\mr^\R(Y_{\HTH}) = 4$, giving
\[\hmrp({\HTH}) = 3 < 4 = \mrp({\HTH}).\]
}\end{ex}

\end{document}